\documentclass{amsart}

\usepackage{amssymb}
\usepackage[all]{xy}
\CompileMatrices

\newtheorem{theorem}{Theorem}[section]
\newtheorem{corollary}[theorem]{Corollary}
\newtheorem{proposition}[theorem]{Proposition}
\newtheorem{lemma}[theorem]{Lemma}

\theoremstyle{definition}
\newtheorem{definition}[theorem]{Definition}
\newtheorem{remark}[theorem]{Remark}

\def\quot{/\!\!/}
\def\mal{\! \cdot \!}
\def\rq#1{\widehat{#1}}
\def\b#1{\overline{#1}}
\def\CC{{\mathbb C}}
\def\TT{{\mathbb T}}
\def\ZZ{{\mathbb Z}}
\def\QQ{{\mathbb Q}}
\def\PP{{\mathbb P}}
\def\SL{{\rm SL}}
\def\Pic{{\rm Pic}}
\def\Supp{{\rm Supp}}
\def\Spec{{\rm Spec}}

\begin{document}

\title[A general Hilbert-Mumford Criterion]
      {A general Hilbert-Mumford Criterion}

\author[J.~Hausen]{J\"urgen Hausen} 
\address{Fachbereich Mathematik und Statistik, Universit\"at Konstanz,
  78457 Konstanz, Germany}
\email{Juergen.Hausen@uni-konstanz.de}
\subjclass{14L24,14L30}

\begin{abstract}
We provide a Hilbert-Mumford Criterion for actions of 
reductive groups $G$ on $\QQ$-factorial complex varieties.
The result allows to construct open subsets admitting 
a good quotient by $G$ from certain maximal open subsets admitting 
a good quotient by a maximal torus of $G$. As an application, we
show how to obtain all invariant open subsets with
good quotient for a given $G$-action on a complete
$\QQ$-factorial toric variety.  
\end{abstract}

\maketitle

\section{Statement of the results}

Let a reductive group $G$ act on a normal complex 
algebraic variety $X$.
It is a central problem in Geometric Invariant Theory 
to construct all $G$-invariant open subsets $V \subset X$ 
admitting a {\em good quotient}, 
i.e.~an affine $G$-invariant morphism $V \to V \quot G$ 
onto a complex algebraic space such that locally
$V \quot G$ is the spectrum of the invariant functions.
Let us call these $V \subset X$ for the moment the 
{\em good $G$-sets}.

In principle, it suffices to know all good $T$-sets
$U \subset X$ for some fixed maximal torus $T \subset G$, 
because the good $G$-sets are precisely the $G$-invariant 
good $T$-sets, see~\cite{BBSw0a}.
The construction of ``maximal'' good $T$-sets is less hard, 
and in order to gain  good $G$-sets 
one studies the following question:
{\em Let $U \subset X$ be a good $T$-set.
When is the intersection $W(U)$ of all translates 
$g \mal U$, $g \in G$, a good $G$-set?}

The classical Hilbert-Mumford Criterion answers this question
in the affirmative for sets of $T$-semistable points of
$G$-linearized ample line bundles. 
Moreover, A.~Bia\l y\-nic\-ki-Birula and J.~\'Swi\c{e}cicka 
settled in~\cite{BBSw0} the case of good $T$-sets defined by 
generalized moment functions, and in~\cite{BBSw0a}
the case $U=X$, as mentioned before.
For $G=\SL_{2}$, several results can be found in~\cite{BBSw1},
\cite{BBSw2}, and~\cite{ha3}. 

As indicated, one imposes maximality conditions on the good 
$T$-set $U$, e.g.~projectivity or completeness of $U \quot T$. 
The most general concept is {\em $T$-maximality}: 
$U$ is not $T$-saturated in some properly larger 
good $T$-set $U'$, 
where {\em $T$-saturated\/} means saturated with respect 
to the quotient map.
For complete $X$ and $T$-maximal $U \subset X$
which are invariant under the normalizer $N(T)$, 
A.~Bia\l ynicki-Birula conjectures that $W(U)$ 
is a good $G$-set~\cite[Conj.~12.1]{BB2}. 

We shall settle the case of {\em $(T,2)$-maximal\/} subsets. 
These are good $T$-sets $U \subset X$ such that 
$U \quot T$ is embeddable into a toric variety, 
and $U$ is not a $T$-saturated subset of some properly larger 
$U'$ having the same properties, compare~\cite{Sw2}.
We shall assume that $X$ is $\QQ$-factorial, 
i.e.~for every Weil divisor on $X$ some 
multiple is Cartier. 
In Section~\ref{sec4}, we prove:

\begin{theorem}\label{maintheorem}
Let a connected reductive group $G$ act on a 
$\QQ$-factorial complex variety~$X$. 
Let $T \subset G$ be a maximal torus and 
$U \subset X$ a $(T,2)$-maximal open subset. 
Then the intersection $W(U)$ of all translates 
$g \mal U$, $g \in G$, is open in $X$, 
there is a good quotient $W(U) \to W(U) \quot G$,
and $W(U)$ is $T$-saturated in $U$. 
\end{theorem}

This generalizes results by 
A.~Bia\l ynicki-Birula and J.~\'Swi\c{e}cicka for 
$X =\PP^{n}$, see~\cite[Thm.~C]{BBSw3},
and by J.~\'Swi\c{e}cicka for smooth complete 
varieties $X$ with $\Pic(X) = \ZZ$, see~\cite[Cor.~6.3]{Sw2}.
As an application of Theorem~\ref{maintheorem}, we obtain:
  
\begin{corollary}\label{toriccase}
Let a connected reductive group $G$ act on a complete
$\QQ$-factorial toric variety $X$, and let $T \subset G$ 
be a maximal torus. Then we have:
\begin{enumerate}
\item For every $T$-maximal open subset $U \subset X$ the set
  $W(U)$ is open and admits a good quotient $W(U) \to W(U) \quot G$. 
\item Every $G$-invariant open subset $V \subset X$ admitting a good 
  quotient $V \to V \quot G$ is a $G$-saturated subset of some set
  $W(U)$ as in~(i). 
\end{enumerate}
\end{corollary}  

Together with well known fan-theoretical descriptions of
the $T$-maximal open subsets, see e.g.~\cite{Sw1}, 
this corollary explicitly solves the quotient problem 
for actions of connected reductive groups $G$ on
$\QQ$-factorial toric varieties.      
In~\cite[Problem~12.9]{BB2} our corollary
was conjectured (in fact for arbitrary toric varieties).

\section{Background on good quotients}\label{sec2}

We recall basic definitions and facts 
on good quotients,
see also~\cite[Chap.~7]{BB2}, \cite[Sec.~1]{BBSw0a} 
and~\cite[Sec.~2]{BBSw3}. 
Let a reductive group $G$ act morphically 
on a complex algebraic variety $X$.
The concept of a good quotient is locally, 
with respect to the \'etale topology,
modelled on the classical invariant theory 
quotient:

\begin{definition}
An $G$-invariant morphism $p \colon X \to Y$ onto
a separated complex algebraic space $Y$ is called a 
{\em good quotient\/} for the $G$-action on $X$ if 
$Y$ is covered by \'etale neighbourhoods $V \to Y$ 
such that
\begin{enumerate}
\item $V$ and its inverse image $U := p^{-1}(V) = X \times_{Y} V$ are 
  affine varieties,
\item $p^{*} \colon \mathcal{O}(V) \to \mathcal{O}(U)$ defines an
  isomorphism onto the algebra of $G$-invariants. 
\end{enumerate}
A good quotient $p \colon X \to Y$ for the $G$-action on $X$
is called {\em geometric}, if its fibres are precisely
the $G$-orbits. 
\end{definition}

A good quotient $X \to Y$ for the $G$-action on $X$ 
is categorical, i.e.~any $G$-invariant morphism $X \to Z$ 
of algebraic spaces factors uniquely through $X \to Y$. 
In particular, good quotient spaces are unique up to 
isomorphism. This justifies the notation
$X \to X \quot G$ for good and $X \to X/G$ 
for geometric quotients.

In the sequel we say that an open subset 
$U \subset X$ of a $G$-variety $X$ with good quotient is
{\em $G$-saturated}, if $U$ is saturated with respect to 
the quotient map $X \to X \quot G$.
The following well known properties of good quotients 
are direct consequences of the corresponding statements 
in the affine case:

\begin{remark}\label{subvarieties}
Assume that the $G$-action on $X$ has a good quotient 
$p \colon X \to X \quot G$. 
\begin{enumerate}
\item If $A \subset X$ is $G$-invariant and closed, then 
  $p(A)$ is closed in $X \quot G$, and the
  restriction $p \colon A \to p(A)$ is a good quotient for the action
  of $G$ on $A$.
\item If $A$ and $A'$ are disjoint $G$-invariant closed subsets of $X$, 
  then $p(A)$ and $p(A')$ are disjoint. 
\item If $U \subset X$ is $G$-saturated and open, then 
  $p(U)$ is open in $X \quot G$, and the restriction 
  $p \colon U \to p(U)$ is a good quotient for the action 
  of $G$ on $U$.
\item If $A \subset X$ and $U \subset X$ are as in~(i) and~(iii), then
  $A \cap U$ is $G$-saturated in $A$.
\end{enumerate}
\end{remark}

Let $X$ be normal (in particular irreducible) 
with a good quotient $X \to X \quot G$. 
Then any reductive subgroup $H \subset G$ admits
a good quotient $X \to X \quot H$, see~\cite[Cor.~10]{BBSw4}.
If $H$ is normal in $G$, then universality of 
good quotients~\cite[Thm.~7.1.4]{BB2}
allows to push down the $G$-action to $X \quot H$. 
Moreover, we have:

\begin{proposition}\label{subgroups}
Let $H \subset G$ be a reductive normal subgroup such that 
$X \quot H$ is an algebraic variety.
Then the canonical map $X \quot H \to X \quot G$ 
is a good quotient for the induced action of $G/H$ on $X \quot H$. 
\endproof
\end{proposition}

We turn to the special case of an action of an algebraic 
torus $T$ on a normal variety $X$. 
Good quotients for such torus actions
are always affine morphisms of normal algebraic varieties,
see~\cite[Cor.~1.3]{BBSw0a}.
We work with the following maximality concepts
for good quotients, compare~\cite[Def.~4.3]{Sw2}: 

\begin{definition}
A $T$-invariant open subset
$U \subset X$ with a good quotient $U \to U \quot T$
is called a {\em $(T,k)$-maximal\/} subset of $X$ if
\begin{enumerate}
\item the quotient space $U \quot T$ is an {\em $A_{k}$-variety}, 
  i.e.~any collection 
  $y_{1}, \ldots, y_{k} \in U \quot T$ admits 
  a common affine neighbourhood in $U \quot T$,
\item $U$ does not occur as proper $T$-saturated subset of some
  $T$-invariant open $U' \subset X$ admitting a good quotient
  $U' \to U' \quot T$ with an $A_{k}$-variety $U' \quot T$.
\end{enumerate}
\end{definition}

As usual, {\em $T$-maximal} stands for $(T,1)$-maximal.
The collection of all $(T,k)$-maximal subsets is always 
finite, see~\cite[Thm.~4.4]{Sw2}. 
The case $k=2$ can also be characterized 
via embeddability of the quotient spaces:
By~\cite[Thm.~A]{Wl}, 
a normal variety has the $A_{2}$-property if and only if
it embeds into a toric variety.

\begin{proposition}\label{toric}
Let $X$ be a toric variety, and let the algebraic torus $T$ act on $X$
via a homomorphism $T \to T_{X}$ to the big torus $T_{X} \subset X$.
Then the $T$-maximal subsets of $X$ are precisely the $(T,2)$-maximal
subsets of $X$.
\end{proposition}

\proof
First observe that every $(T,2)$-maximal subset is 
$T$-saturated in some $T$-maximal subset.
Hence we only have to show that for any $T$-maximal $U \subset X$
the quotient space $U \quot T$ is an $A_{2}$-variety. 
But this is known: 
By~\cite[Cor.~2.4 and~2.5]{Sw1}, the set $U$ is
$T_{X}$-invariant, and $U \quot T$ inherits 
the structure of a toric variety from $U$. 
In particular, $U \quot T$ is an $A_{2}$-variety, 
see~\cite[p.~709]{Wl}.
\endproof

\section{Globally defined $(T,2)$-maximal subsets}\label{sec3}

Let $G$ be a connected reductive group, 
$T \subset G$ a maximal torus, 
and $X$ a normal $G$-variety.
In this section, we reduce the construction of 
$(T,2)$-maximal subsets to a purely 
toric problem in~$\CC^{n}$.  
The following notion is central:

\begin{definition}\label{globaldef}
We say that a $(T,2)$-maximal subset
$U \subset X$ is {\em globally defined\/} in $X$, 
if there are $T$-homogeneous 
$f_{1}, \ldots, f_{r} \in \mathcal{O}(X)$
such that each $X_{f_{i}}$ is an affine open subset of $U$
and any pair $x,x' \in U$ is contained in some $X_{f_{i}}$.
\end{definition}

Here, as usual, $f \in \mathcal{O}(X)$ is called 
{\em $T$-homogeneous}, 
if $f(t \mal x) = \chi(t)f(x)$ holds with 
a character $\chi \colon T \to \CC^{*}$,
and $X_{f}$ denotes the set of all 
$x \in X$ with $f(x) \ne 0$.
Our reduction is split into two lemmas.
The proofs are based on ideas of~\cite{ha2}.

\begin{lemma}\label{reduce2quasiaff}
Let $X$ be $\QQ$-factorial,
and let $U \subset X$ be $(T,2)$-maximal. 
Then there are an algebraic torus $H$ and a 
$\QQ$-factorial quasiaffine $(G \times H)$-variety 
$\rq{X}$ such that
\begin{enumerate}
\item $H$ acts freely on $\rq{X}$ with a $G$-equivariant geometric
  quotient $q \colon \rq{X} \to X$,
\item $\rq{U} := q^{-1}(U)$ is a globally defined
  $(T \times H,2)$-maximal subset of $\rq{X}$.
\end{enumerate}
\end{lemma}

\proof
Let $p \colon U \to U \quot T$ be the quotient.
By assumption, we can cover $U \quot T$ by affine open 
subsets $Y_{1}, \ldots, Y_{r}$ such that any pair 
$y, y' \in U \quot T$ is contained in a common $Y_{i}$. 
Since $p$ is affine, 
each $p^{-1}(Y_{i})$ is affine. Hence each
$X \setminus p^{-1}(Y_{i})$ is of pure codimension one
and, by $\QQ$-factoriality,  
equals the support $\Supp(D_{i})$
of an effective Cartier divisor
$D_{i}$ on $X$.

The Cartier divisors $D_{1}, \ldots, D_{r}$ generate a free 
abelian subgroup $\Lambda$ of the group of all Cartier divisors of 
$X$. Enlarging $\Lambda$ by adding finitely many generators, we 
achieve that every $x \in X$ admits an affine neighbourhood
$X \setminus \Supp(D)$ for some effective member $D \in \Lambda$. 
The group $\Lambda$ gives rise to a
graded $\mathcal{O}_{X}$-algebra:
$$ 
\mathcal{A} 
:= 
\bigoplus_{D \in \Lambda} \mathcal{A}_{D}
:= 
\bigoplus_{D \in \Lambda} \mathcal{O}_{X}(D). 
$$

After eventually replacing $\Lambda$ 
with a subgroup of finite index, 
we can endow $\mathcal{A}$ with a $G$-sheaf structure,
see~\cite[Prop.~3.5]{ha2}:
for any $g \in G$ and any open $V \subset X$,
we then have a $\Lambda$-graded homomorphism  
$\mathcal{A}(V) \to \mathcal{A}(g \mal V)$,
these homomorphisms are compatible with restriction 
of $\mathcal{A}$ and multiplication of $G$,
and the resulting $G$-representation on 
$\mathcal{A}(X)$ is rational. 

We define the desired data; for details see~\cite[Sec.~2]{ha1}.
Let $\rq{X} := \Spec(\mathcal{A})$. The in\-clu\-sion
$\mathcal{O}_{X} \to \mathcal{A}$ 
defines an affine morphism $q \colon \rq{X} \to X$
with $q_{*}(\mathcal{O}_{\rq{X}}) = \mathcal{A}$.
For the canonical section
of an effective $D \in \Lambda$, its zero set in $\rq{X}$
is just $q^{-1}(\Supp(D))$. 
In particular, $\rq{X}$ is covered by affine
sets $\rq{X}_{f}$ and hence is quasiaffine.

The $\Lambda$-grading of $\mathcal{A}$ corresponds to a free
action of the torus $H := \Spec(\CC[\Lambda])$ on $\rq{X}$.
This makes $q \colon \rq{X} \to X$ to an $H$-principal bundle.
In particular, $q$ is a geometric quotient for the $H$-action,
and $\rq{X}$ is $\QQ$-factorial. 
The $G$-sheaf structure of $\mathcal{A}$ induces 
a $G$-action on $\rq{X}$ commuting with the $H$-action
and making $q$ equivariant.

We show that $\rq{U} = q^{-1}(U)$ is $(\rq{T},2)$-maximal,
where we set $\rq{T} := T \times H$. 
First note that the restriction 
$p \circ q \colon \rq{U} \to U \quot T$ 
is a good quotient for the $\rq{T}$-action.
For $(\rq{T},2)$-maximality, let $\rq{U}$ be $\rq{T}$-saturated
in some $(\rq{T},2)$-maximal $\rq{U}_{1} \subset \rq{X}$. 
Then Lemma~\ref{subgroups} gives a commutative diagram
$$
\xymatrix{
{\rq{U}_{1}} \ar[rr]^{{\quot \rq{T}}} \ar[dr]_{/H}^{q}
& &
{\rq{U}_{1} \quot \rq{T}} \\
& {U_{1}} \ar[ur]_{{\quot T}} &}
$$
where $U_{1} := q(\rq{U}_{1})$ is open in $X$. 
Since $\rq{U}$ is $\rq{T}$-saturated in $\rq{U}_{1}$ and 
$\rq{U}_{1} \to U_{1}$ is surjective, this diagram shows
that $U$ is a $T$-saturated subset of $U_{1}$. 
By $(T,2)$-maximality of $U$ in $X$, 
this implies $U=U_{1}$ and hence
$\rq{U} = \rq{U}_{1}$. 

Finally, let $f_{i} \in \mathcal{O}(\rq{X})$ be the canonical 
sections of some large positive multiples of the $D_{i}$. 
The zero set of $f_{i}$ in $\rq{X}$ is just $q^{-1}(\Supp(D_{i}))$. 
In particular, these zero sets are $\rq{T}$-invariant,
and hence the $f_{i}$ are $\rq{T}$-homogeneous. 
By construction, the sets $\rq{X}_{f_{i}}$ equal $q^{-1}(p^{-1}(Y_{i}))$, 
and thus form an affine cover of $\rq{U}$ as required in~\ref{globaldef}. 
\endproof

\begin{lemma}\label{toricextension}
Let $X$ be quasiaffine,
and let $U \subset X$ be a globally defined
$(T,2)$-maximal subset of $X$. 
Then there exist a linear $G$-action on some $\CC^{n}$
and a $G$-equivariant locally closed embedding 
$X \to \CC^{n}$ such that
\begin{enumerate}
\item the maximal torus $T \subset G$ acts on $\CC^{n}$ by means 
  of a homomorphism $T \to \TT^{n}$ to the big torus 
  $\TT^{n} := (\CC^{*})^{n}$,
\item there is a $\TT^{n}$-invariant open $V \subset \CC^{n}$
  containing $U$ as a closed subset and admitting a good quotient
  $V \to V \quot T$.
\end{enumerate}
\end{lemma}

\proof
Let $f_{1}, \ldots, f_{r} \in \mathcal{O}(X)$ be as 
in~\ref{globaldef}, and set $X_{i} := X_{f_{i}}$.
By~\cite[Lemma~2.4]{ha1}, we can realize $X$ as a $G$-invariant
open subset of an affine $G$-variety $\b{X}$
such that the $f_{i}$ extend regularly 
to $\b{X}$ and satisfy $\b{X}_{f_{i}} = X_{i}$. 
Complete the $f_{i}$ to a system 
$f_{1}, \ldots, f_{s}$ of $T$-homogeneous 
generators of the algebra~$\mathcal{O}(\b{X})$.

To proceed, we use the standard representation 
$g \mal f (x) := f(g^{-1} \mal x)$ of $G$
on $\mathcal{O}(\b{X})$. 
Let $M_{i} \subset \mathcal{O}(\b{X})$ be the $G$-module generated by 
$G \mal f_{i}$. Fix a basis $f_{i1}, \ldots, f_{in_{i}}$ of
$M_{i}$ such that all $f_{ij}$ are $T$-homogeneous and for the first
one we have $f_{i1} = f_{i}$. Denoting by $N_{i}$ the dual
$G$-module of $M_{i}$, we obtain $G$-equivariant maps
$$ \Phi_{i} \colon \b{X} \to N_{i}, \qquad x \mapsto [h \mapsto h(x)]. $$

We identify $N_{i}$ with $\CC^{n_{i}}$ by associating to a functional
of $N_{i}$ its coordinates $z_{i1}, \ldots, z_{in_{i}}$ with respect 
to the dual basis $f_{i1}^{*}, \ldots, f_{in_{i}}^{*}$.
Then the pullback $\Phi_{i}^{*}(z_{ij})$ is just the function
$f_{ij}$. Now, consider the direct sum of the $G$-modules
$\CC^{n_{i}}$; we write this direct sum as $\CC^{n}$ but still use
the coordinates $z_{ij}$. 
The maps $\Phi_{i}$ fit together to a
$G$-equivariant closed embedding:
$$ 
\Phi \colon \b{X} \to \CC^{n},
\qquad
x \mapsto 
(f_{11}(x), \ldots, f_{1n_{1}}(x), 
  \ldots, f_{s1}(x), \ldots,  f_{sn_{s}}(x)).
$$

In the sequel, we shall regard $\b{X}$ as a $G$-invariant closed 
subset of $\CC^{n}$. Thus the functions $f_{ij}$ are just the
restrictions of the coordinate functions $z_{ij}$.
By construction, the maximal torus $T$ of $G$ acts
diagonally on $\CC^{n}$, that means that
$T$ acts by a homomorphism $T \to \TT^{n}$
to the big torus $\TT^{n} = (\CC^{*})^{n}$.

We come to the construction of the desired set 
$V \subset \CC^{n}$.
Let $V_{i} \subset \CC^{n}$ be the complement of the 
coordinate hyperplane defined by $z_{i1}$.
Note that $\b{X} \cap V_{i}$ equals $X_{i}$.
In particular, $X_{i}$ is closed in $V_{i}$.
Consider the union $V_{0} := V_{1} \cup \ldots \cup V_{r}$.
Then $V_{0}$ is invariant under the big torus $\TT^{n}$.
Moreover, we have
$$ 
\b{X} \cap V_{0}
\; = \; 
\bigcup_{i=1}^{r} \b{X} \cap V_{i}
\; = \; 
\bigcup_{i=1}^{r} \b{X}_{f_{i}}
\; = \; 
\bigcup_{i=1}^{r} X_{i}
\; = \; 
U.
$$

Let $V \subset V_{0}$ be the minimal $\TT^{n}$-invariant 
open subset with $U = \b{X} \cap V$. Then every closed 
$\TT^{n}$-orbit of $V$ has nontrivial intersection with $U$.
We show that $V$ admits a good quotient by the action of $T$.
By~\cite[Prop.~1.2]{ha2}, it suffices to verify that any 
two points with closed $\TT^{n}$-orbits in $V$ have a 
common $T$-invariant affine open neighbourhood in $V$. 

Let $z,z' \in V$ have closed $\TT^{n}$-orbits in $V$.
Since these $\TT^{n}$-orbits meet $U$,
there are $t,t' \in \TT^{n}$ such that
$t \mal z$ and $t' \mal z'$ lie in $U$. 
By the choice of $f_{1}, \ldots, f_{r}$, 
the points $t \mal z$ and 
$t' \mal z'$ even lie in some common $X_{i}$.
Consider the corresponding $V_{i}$ and the good quotient 
$p \colon V_{i} \to V_{i} \quot T$.
The latter is a toric morphism of affine toric varieties. 

Let $Z_{i} := V_{i} \setminus V$. Then $Z_{i}$ is 
$T$-invariant and closed in $V_{i}$. 
Moreover, $Z_{i}$ does not meet the
$T$-invariant closed subset $X_{i} \subset V_{i}$.
Thus $p(Z_{i})$ and $p(X_{i})$ are closed in 
$V_{i} \quot T$ and disjoint from each other.
In particular, neither $p(t \mal z)$ nor $p(t' \mal z')$ 
lie in $p(Z_{i})$. 
Since $Z_{i}$ is even $\TT^{n}$-invariant, also $p(z)$
and $p(z')$ do not lie in $p(Z_{i})$.

Consequently, there exists a $T$-invariant regular function on
$V_{i}$ that vanishes along $Z_{i}$ but not in the points 
$z$ and $z'$. Removing the zero set of this function from $V_{i}$
yields the desired common $T$-invariant affine open neighbourhood of
the points $z$ and $z'$ in $V$. This proves existence of a good
quotient $V \to V \quot T$.  
\endproof

\section{Proof of the results}\label{sec4}

\proof[Proof Theorem~\ref{maintheorem}] 
First we reduce to the case of globally defined subsets
of quasiaffine varieties. 
So, assume for the moment that Theorem~\ref{maintheorem}
holds in this setting.
Consider the quasiaffine variety $\rq{X}$, 
the torus $H$ and the geometric quotient $q \colon \rq{X} \to X$ 
provided by Lemma~\ref{reduce2quasiaff}.

Then $\rq{G}:= G \times H$ is reductive with maximal torus 
$\rq{T} := T \times H$, and 
$\rq{U} = q^{-1}(U)$ is a globally defined 
$(\rq{T},2)$-maximal subset of~$\rq{X}$. 
By assumption, the intersection $W(\rq{U})$ of all tranlates 
$\rq{g} \mal \rq{U}$ is open, admits a good 
quotient by $\rq{G}$, and is 
$\rq{T}$-saturated in $\rq{U}$.
Since each $\rq{g} \mal \rq{U}$ is $H$-invariant 
and $q \colon \rq{X} \to X$ is $G$-equivariant, we obtain
$$ 
W(\rq{U})
\; = \; 
\bigcap_{\rq{g} \in \rq{G}} \rq{g} \mal \rq{U}
\; = \; 
\bigcap_{g \in G} g \mal \rq{U}
\; = \; 
\bigcap_{g \in G} g \mal q^{-1}(U)
\; = \; 
q^{-1}(W(U)).
$$ 

In particular, $W(U)$ is open in $X$.
Moreover, restricting $q$ gives a geometric quotient 
$W(\rq{U}) \to W(U)$ for the $H$-action.
Lemma~\ref{subgroups} tells us that the
induced map from $W(U)$ onto $W(\rq{U}) \quot \rq{G}$ 
is a good quotient for the $G$-action on 
$W(U)$.
Similarly, we infer $T$-saturatedness of $W(U)$ in $U$
from the commutative diagram
$$
\xymatrix{
{\rq{U}} \ar[rr]^{{\quot \rq{T}}} \ar[dr]_{/H}^{q}
& &
{\rq{U}} \quot \rq{T} \\
& {U} \ar[ur]_{{\quot T}} &}
$$

We are left with proving~\ref{maintheorem} for quasiaffine $X$
and globally defined $(T,2)$-maximal $U \subset X$.
By Lemma~\ref{toricextension}, we may view 
$X$ as a $G$-invariant locally closed subset of
a $G$-module $\CC^{n}$, where
$T$ acts via a homomorphism $T \to \TT^{n}$
and $U$ is closed in some $\TT^{n}$-invariant 
open $V \subset \CC^{n}$
with good quotient $V \to V \quot T$.
We regard $\CC^{n}$ as the $G$-invariant open subset 
of $\PP^{n}$ obtained by removing the zero
set of the homogeneous coordinate $z_{0}$. 

Let $V' \subset \PP^{n}$ be a $T$-maximal open subset
containing $V$ as a $T$-saturated subset. 
Let $\b{X}$ be the closure of $X$ in $\PP^{n}$, 
and set $X' := \b{X} \cap V'$. 
Then $X'$ is closed in $V'$, and we have $U = X' \cap V$.
Using~\ref{subvarieties}~(i), (iii) and~(iv),
we subsume the situation in a commutative cube
$$
\xymatrix@!0{
& 
U \ar[rr]  \ar[dl] \ar'[d][dd]
& & 
V \ar[dd] \ar[dl]
\\
X' \ar[rr] \ar[dd]
& & 
V'  \ar[dd]
\\
&
{U \quot T} \ar[dl] \ar'[r][rr]
& & 
{V \quot T} \ar[dl]
\\
{X' \quot T} \ar[rr] 
& & 
{V' \quot T} 
}
$$
where the downwards arrows are good quotients 
by the respective actions of $T$, 
the right arrows are closed inclusions,
the upper diagonal arrows are $T$-saturated inclusions
and the lower diagonal arrows are open inclusions.

According to~\cite[Thm.~C]{BBSw3}, 
the intersection $W(V')$ of all translates $g \mal V'$
is open in $\PP^{n}$ and admits a good quotient by the 
action of $G$. 
Recall from~\cite[Lemma~8.4]{BBSw3} that $W(V')$ is 
$T$-saturated  in $V'$.
We transfer the desired properties step by step from 
$W(V')$ to $W(U)$. 
First note that by $G$-invariance of $\b{X}$ we have
$$
W(X') 
\; = \;
\bigcap_{g \in G} g \mal X'
\; = \;
\bigcap_{g \in G} g \mal (\b{X} \cap V')
\; = \;
\b{X} \cap W(V')
\; = \;
X' \cap W(V').
$$ 

Thus $W(X')$ is open in $X'$, and by~\ref{subvarieties}~(iv) it is 
$T$-saturated in $X'$.
In particular, the $T$-action on $W(X')$ has a good quotient.
Moreover, $W(X')$ is $G$-invariant and closed in $W(V')$.
Thus~\ref{subvarieties}~(i)
ensures the existence of a good quotient 
$$u \colon W(X') \to W(X') \quot G.$$

Consider $B := X' \setminus X$. Since $X$ is open 
in $\b{X}$ and $B$ equals $(\b{X} \setminus X) \cap X'$, 
the set $B$ is closed in~$X'$. 
The intersection $W(B)$ of the translates $g \mal B$, 
where $g \in G$, is $G$-invariant and closed in $W(X')$. 
We claim that it suffices to verify
\begin{equation}\label{eqn1}
W(U) = W(X') \setminus u^{-1}(u(W(B))).
\end{equation}

Indeed, suppose we have (\ref{eqn1}). Then $W(U)$ 
is open in $X'$, hence in~$U$, and thus in~$X$. 
Property~\ref{subvarieties}~(iii) provides
a good quotient $W(U) \to W(U) \quot G$. 
Moreover, $W(U)$ is $T$-saturated in $W(X')$,
because it is $G$-saturated and we have 
the induced map from $W(X') \quot T$ onto $W(X') \quot G$.
Since $W(X')$ and $U$ are $T$-saturated in $X'$, 
we obtain that $W(U)$ is $T$-saturated in $U$.

We verify~(\ref{eqn1}). Let 
$v \colon X' \to X' \quot T$ be the quotient map.
As a subvariety, $X' \quot T$ inherits the
$A_{2}$-property from $V' \quot T$, which in turn 
satisfies it by~\ref{toric}. 
Thus, since $U$ is $(T,2)$-maximal in $X$, it is necessarily
the maximal $T$-saturated subset of $X'$ which is contained 
in $X \cap X'$. 
In terms of $B = X' \setminus X$ this means:
\begin{equation}\label{eqn2}
U = X' \setminus v^{-1}(v(B)). 
\end{equation}

We check the inclusion ``$\subset$'' of~(\ref{eqn1}).
Let $x \in u^{-1}(u(W(B)))$.
Then, by~\ref{subvarieties}~(ii), the 
closure of $G \mal x$ meets $W(B)$.
The classical Hilbert-Mumford Lemma~\cite[Thm.~4.2]{Bi} 
says that for some maximal torus $T' \subset G$
the closure of $T' \mal x$ meets $W(B)$. 
Let $g \in G$ with $T = gT'g^{-1}$. Then the closure of 
$T \mal g \mal x$ meets $W(B)$. 
Hence $g \mal x$ lies in $v^{-1}(v(B))$.
By~(\ref{eqn2}), the point $x$ cannot belong to~$W(U)$.

We turn to the inclusion ``$\supset$'' of~(\ref{eqn1}).
For this, consider the set 
$A := (X \cap X') \setminus U$. 
Then $X'$ is the disjoint union of $U$, $A$ and $B$.
Consequently, we have
$$
W(U)
\; = \; 
\bigcap_{g \in G} g \mal (X' \setminus (A \cup B))
\; = \; 
W(X') \setminus 
\bigcup_{g \in G} g \mal A \cup g \mal B.
$$

So we have to show that $u$ maps a given
$x \in W(X') \cap g \mal (A \cup B)$ to $u(W(B))$. 
Since $g^{-1} \mal x \not\in U$ holds, we infer from~(\ref{eqn2})
that $g^{-1} \mal x$ lies in $v^{-1}(v(B))$. 
According to~\ref{subvarieties}~(ii), the
closure of $T \mal g^{-1} \mal x$ in $X'$ meets $B$. Since 
$W(X')$ is $T$-saturated in $X'$, this implies that the closure of 
$T \mal g^{-1} \mal x$ meets $W(X') \cap B$. But we have
$$ 
W(X') \cap B
\; = \; 
W(X') \setminus X
\; = \; 
\bigcap_{g \in G} g \mal (X' \setminus X)
\; = \; 
W(B).
$$

Hence we obtained that the closure of the orbit $G \mal x$ 
intersects $W(B)$. This in turn shows that the image $u(x)$ 
lies in $u(W(B))$. 
\endproof

\proof[Proof Corollary~\ref{toriccase}]
Recall from~\cite[Sec.~4]{Co} that the automorphism group of 
$X$ is a linear algebraic group having the big torus 
$T_{X} \subset X$ as a maximal torus. 
Thus, by conjugating $T_{X}$ we achieve that  
$T \subset G$ acts on $X$ 
via a homomorphism $T \to T_{X}$. 
Proposition~\ref{toric} then ensures that each $T$-maximal 
subset of $X$ is as well $(T,2)$-maximal, 
and statement~(i) follows from Theorem~\ref{maintheorem}.

For statement~(ii), let $V \subset X$ be open and 
$G$-invariant with good quotient $V \to V \quot G$. 
Then~\cite[Cor.~10]{BBSw4} provides a good quotient 
$V \to V \quot T$. 
Let $U \subset X$ be a $T$-maximal subset containing 
$V$ as $T$-saturated subset.
Then we have $V \subset W(U)$.
Again by~\ref{toric}, the set $U$ is $(T,2)$-maximal.
Thus Theorem~\ref{maintheorem} says that $W(U)$ is open,
has a good quotient $u \colon W(U) \to W(U) \quot G$, and is
$T$-saturated in $U$.

For $G$-saturatedness of $V$ in $W(U)$ 
we have to show that any $x \in u^{-1}(u(V))$ with
closed $G$-orbit in $W(U)$ belongs to $V$. 
For this note that $V$ is
$T$-saturated in $W(U)$, because both sets are so in $U$. 
Now, let $y \in V$ with $u(y) = u(x)$. 
Then $x$ lies in the closure of $G \mal y$.
Thus~\cite[Thm.~4.2]{Bi}
provides a $g \in G$ such that
the closure of $T \mal g \mal y$ meets $G \mal x$. 
Since $g \mal y$ lies in $V$ and $V$ is $T$-saturated in 
$W(U)$, we obtain $G \mal x \subset V$, and hence $x \in V$.
\endproof

\bibliography{}

\begin{thebibliography}{}%
%
\bibitem{BB2} A.~Bia\l ynicki-Birula: Algebraic
  Quotients. In: R.V.~Gamkrelidze, V.L.~Popov (Eds.), Encyclopedia of
  Mathematical Sciences, Vol. 131., 1--82 (2002) 
%
\bibitem{BBSw0}  A.~Bia\l ynicki-Birula, J.~\'Swi\c{e}cicka:
  Generalized moment functions and orbit spaces. Amer. J. Math.,
  Vol.~109, 229--238 (1987)
%
\bibitem{BBSw0a}  A.~Bia\l ynicki-Birula, J.~\'Swi\c{e}cicka: A
  reduction theorem for existence of good quotients. Amer.~J.~Math.,
  Vol.~113, 189--201 (1990)
%
\bibitem{BBSw1}  A.~Bia\l ynicki-Birula, J.~\'Swi\c{e}cicka: On
  complete orbit spaces of $\SL(2)$-actions. Colloq. Math. Vol.~55,
  No. 2, 229--241 (1988)
%
\bibitem{BBSw2}  A.~Bia\l ynicki-Birula, J.~\'Swi\c{e}cicka: On
  complete orbit spaces of $\SL(2)$-actions II. 
  Colloq. Math. Vol.~63, No.~1, 9--20 (1992)
%
\bibitem{BBSw3}  A.~Bia\l ynicki-Birula, J.~\'Swi\c{e}cicka: Open
  subsets of projective spaces with a good quotient by an action of a
  reductive group. Transform. Groups, Vol.~1, No.~3, 153--185 (1996).
%
\bibitem{BBSw4}  A.~Bia\l ynicki-Birula, J.~\'Swi\c{e}cicka: Three
  theorems on existence of good quotients. Math. Ann., Vol.~307,
  143--149 (1997).
%
\bibitem{Bi} D.~Birkes: Orbits of linear algebraic groups. 
  Ann. Math., Ser~2, Vol.~93, 459--475 (1971).
%
\bibitem{Co} D. Cox: The homogeneous coordinate ring of a toric
  variety. J. Algebr. Geom., Vol. 4, No. 1, 17--50 (1995). 
%
\bibitem{ha1} J. Hausen: Equivariant embeddings into smooth toric
  varieties. Canad. Math. J., Vol.~54, No.~3, 554--570 (2002)
%
\bibitem{ha2} J. Hausen: Producing good quotients by embedding into 
  toric varieties. S\'emin. et Congr\`es~6 (SMF), 193--212 (2002)  
%
\bibitem{ha3} J. Hausen: A Hilbert-Mumford Criterion for
  $\rm{SL}_{2}$-actions. Preprint, Konstanzer Schriften in Mathematik
  und Informatik~166 (2002)
%
\bibitem{Sw1} J.~\'Swi\c{e}cicka: Quotients of toric varieties by
  actions of subtori. Colloq. Math., Vol.~82, No.~1, 105--116 (1999). 
%
\bibitem{Sw2} J.~\'Swi\c{e}cicka: A combinatorial construction of sets
  with good quotients by an action of a reductive
  group. Colloq. Math., Vol.~87, No.~1, 85--102 (2001).
%
\bibitem{Wl} J.~W\l odarczyk: Embeddings in toric varieties and
  prevarieties. J. Algebr. Geom., Vol.~2, No.~4, 705--726 (1993).
\end{thebibliography}

\end{document}